\documentclass[11pt,a4paper]{article}
\usepackage[T1]{fontenc}
\usepackage[utf8]{inputenc}
\usepackage[english]{babel}
\usepackage{lmodern}
\usepackage{microtype}
\usepackage[margin=27mm,headheight=14pt]{geometry}
\usepackage{amsmath,amssymb,amsthm,mathtools}
\usepackage{tikz,tikz-cd}
\usetikzlibrary{arrows.meta,calc}
\usepackage{enumitem}
\usepackage{xcolor}
\usepackage{needspace}
\usepackage{hyperref}
\hypersetup{colorlinks=true,linkcolor=blue!45!black,citecolor=blue!45!black,urlcolor=blue!45!black,pdftitle={Hidden Symmetries of the Kepler Problem and Jacobi's Theorem},pdfauthor={Juan B. Sancho}}
\newtheorem{thm}{Theorem}[section]
\newtheorem{theorem}[thm]{Theorem}
\newtheorem{proposition}[thm]{Proposition}
\newtheorem{corollary}[thm]{Corollary}

\theoremstyle{definition}
\newtheorem{observation}[thm]{Observation}

\newtheorem*{definition}{Definition}
\theoremstyle{remark}
\newtheorem*{remark}{Remark}

\setlist[enumerate]{itemsep=5pt,topsep=5pt}
\allowdisplaybreaks[2]
\title{Hidden Symmetries of the Kepler Problem\\and Jacobi's Theorem}
\author{ Alejandro Moreno D\'{i}az\\ {\small {\it E-mail:} amorenop@unex.es}
 \and 
 Juan B. Sancho de Salas\\
 {\small {\it E-mail:} jsancho@unex.es}\smallskip
 \\{\small Dpto. de Matem\'{a}ticas, Univ. de Extremadura, E-06071 Badajoz (SPAIN)}}
\begin{document}
\maketitle

\begin{abstract} We present a geometric derivation of the hidden symmetries of the Kepler problem, combining a duality for central-force systems, obtained by interchanging positions and velocities, with a generalized form of Jacobi's theorem. Applied to the dual Kepler problem, this theorem identifies Kepler hodographs of prescribed energy with geodesics of a metric of constant curvature. After adjoining one point at infinity, the resulting geometry is complete. Its Killing vector fields induce a finite-dimensional Lie algebra of infinitesimal symmetries on the corresponding energy hypersurface. The distinguished point at infinity determines a canonical decomposition of the Killing algebra into two vector subspaces, whose associated Noether invariants recover angular momentum and the Runge-Lenz vector. 
\end{abstract}

\section*{Introduction}

The Kepler problem (the Newtonian two-body problem) in $\mathbb R^3$ has some rather mysterious properties. It possesses ``hidden'' symmetries, which do not arise from the geometry of configuration space and are associated with the so-called Runge--Lenz vector. In the literature, this vector is generally introduced unexpectedly, without a clear motivation (see, for example, \cite{Goldstein,Thirring2003}). The angular momentum and Runge--Lenz vectors generate a Lie algebra of symmetries which, by a direct calculation, is isomorphic to $\mathfrak{so}(4)$ on each negative energy level. The purpose of this paper is to show that the surprising appearance of $\mathfrak{so}(4)$ can be explained by a classical theorem of Jacobi.

Jacobi's theorem states that the trajectories of a conservative mechanical system with a prescribed energy $e$ coincide, up to reparametrization, with the geodesics of a Riemannian metric. Moreover, the Killing vector fields of this Riemannian metric define symmetries of the mechanical system and hence Noether invariants of the motion at the prescribed energy.

We shall show that the hidden symmetries of the Kepler problem can be obtained naturally by applying Jacobi's theorem, not to the Kepler problem itself, but to what we shall call the ``dual problem.'' We briefly outline the procedure.

\begin{enumerate}
\item In $n$-dimensional Euclidean space $(\mathbb R^n,g=\sum_i \text{d}x_i^2)$, consider a central-force mechanical system whose trajectories satisfy a differential equation of the form
\begin{equation}\label{eq:central-intro}
\mathbf x''=h(\mathbf x,\mathbf x')\mathbf x.
\end{equation}
For the Kepler problem, the differential equation is
$$
\mathbf x''=-k\frac{\mathbf x}{|\mathbf x|^3},
$$
where $k$ is a positive constant. If we apply the anti-symplectic diffeomorphism
$$
J:T\mathbb R^n\longrightarrow T\mathbb R^n,
\qquad (\mathbf x,\mathbf v)\longmapsto(\mathbf v,\mathbf x),
$$
which interchanges positions and velocities, to the trajectories of equation~\eqref{eq:central-intro} in phase space $T\mathbb R^n$, we obtain the reparametrized trajectories of another differential equation,
$$
\mathbf x''=\frac{1}{h(\mathbf x',\mathbf x)}\mathbf x,
$$
which we call the dual problem. For the Kepler equation, the differential equation of the dual problem is
$$
\mathbf x''=-\frac{|\mathbf x'|^3}{k}\mathbf x.
$$
Geometrically, the trajectories of the dual problem are the reparametrized hodographs of the original problem. Furthermore, the symmetries of the two problems, and their corresponding Noether invariants, are related by the diffeomorphism $J$.

\item Applying Jacobi's theorem to the dual problem, we find that its trajectories, for fixed mass $m$ and energy $e$, coincide with the reparametrized geodesics of the Riemannian metric
\[
\bar g_e:=\frac{4}{(|\mathbf x|^2+K)^2}\sum_i \text{d}x_i^2,
\qquad K=-\frac{2e}{m},
\]
on the region $|\mathbf x|^2+K>0$ in $\mathbb R^n$.

\item The Riemannian metric $\bar g_e$ has constant curvature $K$. After adjoining a point $p_\infty$ at infinity in $\mathbb R^n$, this geometry is isometric, through a suitable stereographic projection, to one of the three classical geometries $\mathbb{S}^n$, $\mathbb R^n$, and $\mathbb{H}^n$, according to the sign of $K$. Thus its Lie algebra of Killing vector fields is $\mathfrak{so}(n+1)$ for $K>0$ (that is, $e<0$), $\mathfrak{so}(1,n)$ for $K<0$ (that is, $e>0$), or $\mathfrak{se}(n)$ for $K=0$ (that is, $e=0$). In this way we obtain a finite-dimensional Lie algebra of the Kepler problem on the hypersurface of  energy $e$.

\item The Lie algebra of Killing vector fields of $\bar g_e$ splits canonically as a direct sum of two vector subspaces. The first is the Lie subalgebra of Killing vector fields vanishing at $p_\infty$. The second is the vector subspace of Killing vector fields $K$ satisfying $(\nabla K)_{p_\infty}=0$. Through the corresponding Noether invariants, this second subspace determines the Runge--Lenz vector, which therefore arises intrinsically.
\end{enumerate}

Historically, the study of Kepler hodographs began with Hamilton's discovery in 1846 that these curves are circular arcs or intervals of straight lines through the origin. In the negative-energy case $e<0$, Moser~\cite{Moser1970} discovered that the hodographs are stereographic projections of great circles on a sphere $\mathbb{S}^n$; see also \cite{Gyorgyi1968,Milnor1983}. Moser's result was extended to positive and zero energies by Osipov~\cite{Osipov1972,Osipov1977} and Belbruno~\cite{Belbruno1977}, with the sphere $\mathbb{S}^n$ replaced by hyperbolic space $\mathbb{H}^n$ and Euclidean space $\mathbb R^n$. Osipov~\cite{Osipov1977} also established the relationship between the symmetries of these constant-curvature geometries and the symmetries of the Kepler problem. The metric $\bar g_e$ appears in Osipov~\cite{Osipov1977} for $m=1$ (and in Moser~\cite{Moser1970} for $e<0$); see also \cite{Chanda2016} and \cite{Marle2012}. The novel aspect here is the naturalness with which the metric $\bar g_e$ is obtained, using a SODE duality and Jacobi's theorem.

The classical theorem of Jacobi applies to conservative mechanical systems. This might appear to present a difficulty, since the dual of the Kepler problem is not conservative. Fortunately, we shall prove that Jacobi's theorem and the theory of Noether invariants remain valid for more general mechanical systems, which we call \emph{quasi-conservative}; see Sections~2 and~3. In these systems, the dynamical vector field on phase space is proportional to a Hamiltonian vector field.

The paper is organized as follows. Section 1 reviews the geometric theory of classical mechanics. Sections 2 and 3 reformulate Noether invariants and the Jacobi theorem for quasi-conservative systems. Section 4 examines a notion of duality for mechanical systems with central forces. Finally, we apply the previous elements to the Kepler problem, obtaining the hidden symmetries and the Runge-Lenz vector.

\section{Preliminaries}

To make this article self-contained and to establish the notation, in this section we review the geometric theory of autonomous classical mechanics.

\subsection*{The symplectic structure of the tangent bundle}

Throughout this section, $(Q,g)$ is a Riemannian manifold and $\pi\colon TQ\to Q$ is its tangent bundle. Physicists use the symbol $Q$ for the configuration space of a mechanical system.

\paragraph{Coordinates on $TQ$.}
Let $(V,q_1,\ldots,q_n)$ be a coordinate chart on $Q$; physicists call the $q_i$ generalized coordinates. We define coordinates $q_1,\ldots,q_n,v_1,\ldots,v_n$ on the open set $TV\subseteq TQ$ by
\[
q_i(v_q):=q_i(q),\qquad v_i(v_q):=v_q(q_i).
\]
Thus the $2n$ coordinates of a vector $v_q\in TV$ consist of the $n$ coordinates of its base point $q$, together with the $n$ components of $v_q$ in the basis $(\partial_{q_1},\ldots,\partial_{q_n})_q$.

\begin{definition}
The \emph{universal vector field} is the vector field $\mathbb{U}$ along $\pi$, that is, the section of $\pi^*TQ\to TQ$ defined by
$$
\mathbb{U}_{v_q}:=\,v_q\qquad\text{for every }v_q\in TQ.
$$
\end{definition}
In the local coordinates $(q_i,v_i)$ defined above,
$$
\mathbb{U}\,=\,\sum_i v_i\partial_{q_i}.
$$
The cotangent bundle of a manifold $Q$ carries a canonical $1$-form, the Liouville form, and a canonical symplectic structure. If $Q$ is equipped with a Riemannian metric, these forms can be transported to the tangent bundle by the metric isomorphism $TQ\simeq T^*Q$. The resulting forms on $TQ$ can also be defined directly as follows.

\begin{definition}
The \emph{Liouville form} is the $1$-form $\Theta$ on $TQ$ defined by $\,\Theta(D)=\mathbb{U}\cdot\pi_*D$, where the scalar product is taken with respect to $g$. Explicitly,
\[
\Theta_{v_q}(D_{v_q})
:=\,\mathbb{U}_{v_q}\cdot\pi_*(D_{v_q})
\,=\,v_q\cdot\pi_*(D_{v_q}).
\]
\end{definition}
In the coordinates $(q_i,v_i)$,
\[
\Theta\,=\,\sum_{i,j}g_{ij}v_j\,\text{d}q_i\,=\,\sum_i p_i\,\text{d}q_i,
\]
where the functions $p_i:=\sum_jg_{ij}v_j$ are called \emph{generalized momenta}.

\begin{definition}
The \emph{canonical symplectic form} on $TQ$ is the $2$-form
\[
\omega:=\,\text{d}\Theta.
\]
\end{definition}
In the coordinates $(q_i,p_i)$,
\[
\omega\,=\,\sum_i \text{d}p_i\wedge \text{d}q_i,
\]
so $\omega$ is indeed a symplectic form on $TQ$.

\subsection*{Second-order differential equations}

Every curve $\sigma\colon I\to Q$ has a natural lift to the tangent bundle,
\[
\begin{tikzcd}
& TQ\arrow[d,"\pi"]\\
I\arrow[ur,"\widetilde\sigma"]\arrow[r,"\sigma"']&Q,
\end{tikzcd}
\]
defined by $\widetilde\sigma(t)=v_t$, where $v=\sigma_*\partial_t$ is the tangent vector field of $\sigma$.

In local coordinates, the lift of $\sigma(t)=(q_1(t),\ldots,q_n(t))$ is
\[
\widetilde\sigma(t)
=\bigl(q_1(t),\ldots,q_n(t),v_1(t)=\dot q_1(t),\ldots,v_n(t)=\dot q_n(t)\bigr).
\]

\begin{definition}
A \emph{SODE} (second-order differential equation) on a manifold $Q$ is a vector field $Z$ on $TQ$ satisfying $\pi_*Z=\mathbb{U}$; equivalently,
\[
\pi_*(Z_{v_q})\,=\,v_q\, .
\]
\end{definition}
In local coordinates, this property means that $Z$ has the form
\[
Z\,=\,\sum_i v_i\partial_{q_i}+\sum_i f_i(q,v)\partial_{v_i}.
\]
To justify the definition, observe that finding the integral curves $(q(t),v(t))$ of $Z$ amounts to solving the first-order system
\[
\begin{cases}
\dot q_i(t)\,=\,v_i(t),\\
\dot v_i(t)\,=\,f_i(q(t),v(t)),
\end{cases}
\]
which is equivalent to the second-order system
\[
\ddot q_i(t)\,=\,f_i(q(t),\dot q(t)).
\]

\begin{definition}
A \emph{solution} of a SODE $Z$ (a vector field on $TQ$) is a curve $\sigma:I\to Q$ whose lift $\widetilde\sigma:I\to TQ$ is an integral curve of $Z$.
\end{definition}

\subsection*{Autonomous mechanical systems}

\begin{definition}
An \emph{autonomous mechanical system} is a triple $(Q,g,F)$, where $Q$ is a differentiable manifold, called the \emph{configuration space}; $g$ is a Riemannian metric on $Q$, called the \emph{kinetic metric}; and $F$ is a vector field along $\pi:TQ\to Q$, called the \emph{force-intensity field}.
\end{definition}
The \emph{kinetic energy} of a mechanical system $(Q,g,F)$ is the function on its tangent bundle given by
\[
\mathsf{T}:TQ\longrightarrow\mathbb R,
\qquad
\mathsf{T}(v_q):=\,\tfrac12g_q(v_q,v_q)\,=\,\tfrac12|v_q|^2.
\]
The kinetic metric and the kinetic energy encode the same information, up to the factor $1/2$.

\begin{definition}
The tangent bundle $TQ$ of a mechanical system $(Q,g,F)$ is called its \emph{phase space}. It carries the Liouville $1$-form $\Theta$ and the symplectic $2$-form $\omega=\text{d}\Theta$ associated with the Riemannian structure of $(Q,g)$.
\end{definition}
Let $\nabla$ denote the Levi-Civita connection of $(Q,g)$.

\begin{definition}
The \emph{solutions} of a mechanical system $(Q,g,F)$ are the curves $\sigma\colon I\to Q$ whose tangent vector field $v=\sigma_*\partial_t$ satisfies
\[
(\nabla_{\partial_t}v)_t\,=\,F_{v_t}\ .
\]
\end{definition}
In particular, in the absence of forces, $F=0$, the solutions are the geodesics of the configuration space $(Q,g)$.

\begin{definition}
For a mechanical system $(Q,g,F)$, the \emph{dynamical vector field} is the unique vector field $Z$ on phase space $TQ$ whose integral curves $\bar\sigma\colon I\to TQ$ are the lifts of the solution curves $\sigma\colon I\to Q$ of the mechanical system.
\end{definition}
Equivalently, $Z$ is the unique vector field on $TQ$ satisfying
\[
\pi_*Z\,=\,\mathbb{U},\qquad \nabla_Z\mathbb{U}\,=\,F\ ,
\]
as equalities of vector fields along $\pi\colon TQ\to Q$.

The condition $\pi_*Z=\mathbb{U}$ says that the dynamical vector field is a SODE. Conversely, every SODE on a Riemannian manifold $(Q,g)$ can be regarded as the dynamical vector field of a mechanical system by defining the force intensity through the equation $\nabla_Z\mathbb{U}=F$.

\begin{definition}
The \emph{work form} of a mechanical system $(Q,g,F)$ is the $1$-form $\textsc{w}$ on $TQ$ defined by
\[
\textsc{w}_{v_q}(D_{v_q}):=\,F_{v_q}\cdot\pi_*(D_{v_q})\ ,
\]
where the scalar product is taken with respect to the kinetic metric $g$.
\end{definition}

\begin{theorem}\label{thm:mechanical-field}
The dynamical vector field of a mechanical system $(Q,g,F)$ is the unique vector field $Z$ on $TQ$ satisfying
$$
i_Z\omega\,=\,\textsc{w}-\mathrm{d}\mathsf{T}.
$$
\end{theorem}
\begin{proof}
Suppose first that $Q$ is an open subset of Euclidean space, with its standard metric. In the induced coordinates $(q^i,v^i)$ on $TQ$, the field of motion and the symplectic form are
$$
Z\,=\,\sum_i\left(v^i\frac{\partial}{\partial q^i}
+F^i\frac{\partial}{\partial v^i}\right),
\qquad
\omega\,=\,\sum_i \text{d}v^i\wedge \text{d}q^i.
$$
Since
$$
\textsc{w}\,=\,\sum_i F^i \text{d}q^i,
\qquad
\mathsf{T}\,=\,\frac12\sum_i(v^i)^2,
$$
direct contraction gives
$$
i_Z\omega\,
=\,\sum_i\left(F^i\text{d}q^i-v^i\text{d}v^i\right)\,
=\,\textsc{w}-\text{d}\mathsf{T}.
$$

Now let $(Q,g)$ be an arbitrary Riemannian manifold. Fix $q\in Q$ and choose normal coordinates centered at $q$. Then
$$
g_{ij}(q)\,=\,\delta_{ij}\ ,
\qquad
\frac{\partial g_{ij}}{\partial q^k}(q)=0\ .
$$
Thus the metric agrees with the Euclidean metric to first order at $q$, and its Christoffel symbols vanish there.

Consequently, at every $v_q\in T_qQ$, the coordinate expression of $Z_{v_q}$ is the same as in the Euclidean case. Moreover,
$$
\omega\,=\,\text{d}\left(\sum_{i,j}g_{ij}v^j\text{d}q^i\right)\ ,
\qquad
\mathsf{T}\,=\,\frac12\sum_{i,j}g_{ij}v^iv^j\ ,
\qquad
\textsc{w}\,=\,\sum_{i,j}g_{ij}F^j\text{d}q^i\ .
$$
The values of $\omega$, $\text{d}\mathsf{T}$, and $\textsc{w}$ at $v_q$ therefore also coincide with their Euclidean expressions. The preceding computation proves
$$
(i_Z\omega)_{v_q}\,=\,(\textsc{w}-\text{d}\mathsf{T})_{v_q}.
$$
Since $q$ and $v_q$ were arbitrary, the identity holds throughout $TQ$.

Finally, the nondegeneracy of $\omega$ implies that a vector field satisfying this identity is unique.
 
\end{proof}

\begin{definition}
A mechanical system $(Q,g,F)$ is \emph{conservative} if its work form is exact: $\textsc{w}=-\text{d}\,\mathsf{U}$.
\end{definition}
The function $\mathsf{U}$ on $TQ$ is called the \emph{potential energy}. The \emph{total energy} is the sum $\textsf{E}=\mathsf{T}+\mathsf{U}$ of the kinetic and potential energies.

\begin{remark}
In a conservative system $(Q,g,F)$, the potential energy $\mathsf{U}$ is a function on $Q$. Consequently, the force intensity $F=-\operatorname{grad}\mathsf{U}$ is also defined on $Q$.

Although $\mathsf{U}$ has initially been defined as a function on $TQ$, it is in fact a function on $Q$; more precisely, it is the pullback of a function on $Q$ under $\pi\colon TQ\to Q$. For every vertical vector field $D$ on $TQ$, meaning that $D$ is tangent to the fibers of $\pi$ or, equivalently, that $\pi_*D=0$, we have
\[
D(\mathsf{U})\,=\,(\text{d}\,\mathsf{U})(D)\,=\,-\textsc{w}(D)\,
=\,-F\cdot\pi_*D\,=\,F\cdot0\,=\,0 .
\]
Thus $\mathsf{U}$ is constant along the fibers of $\pi\colon TQ\to Q$ and is therefore the pullback of a function on $Q$.
\end{remark}

\begin{corollary}\label{cor:hamilton-equations}
The dynamical vector field of a conservative mechanical system is the unique vector field $Z$ on $TQ$ satisfying \textbf{Hamilton's equations}:
\[
i_Z\omega\,=\,-\mathrm{d}\textsf{E}.
\]
\end{corollary}
\begin{proof}
$i_Z\omega=\textsc{w}-\text{d}\mathsf{T}=-\text{d}\,\mathsf{U}-\text{d}\mathsf{T}=-\text{d}E$.

\end{proof}

Consider the coordinates $(q_1,\ldots,q_n,p_1,\ldots,p_n)$ on $TQ$ in which the symplectic form is $\omega=\sum_i \text{d}p_i\wedge\text{d}q_i$. Writing the trajectories $\bar\sigma(t)=(q(t),p(t))$ of the dynamical vector field in these coordinates, so that $Z=\sum_i\dot q_i\partial_{q_i}+\sum_i\dot p_i\partial_{p_i}$ along a trajectory, Hamilton's equations take their classical form:
\[
\dot q_i\,=\,\frac{\partial \textsf{E}}{\partial p_i}\ ,\qquad
\dot p_i\,=\,-\frac{\partial \textsf{E}}{\partial q_i}\ .
\]

\begin{corollary}\label{cor:energy-conservation}
In a conservative mechanical system, the total energy is constant along the trajectories of the motion: $Z(\mathsf{E})=0$.
\end{corollary}
\begin{proof}
By Corollary~\ref{cor:hamilton-equations},
$$
Z(\textsf{E})\,=\,(\text{d}\mathsf{E})(Z)\,=\,-(i_Z\omega)(Z)\,=\,-\omega(Z,Z)\,=\,0\ .
$$
\end{proof}

\begin{corollary}\label{cor:lie-work}
The dynamical vector field $Z$ of a mechanical system $(Q,g,F)$ satisfies
\[
L_Z\omega\,=\,\mathrm{d}\textsc{w}\ .
\]
\end{corollary}
\begin{proof}
Cartan's formula gives
\[
L_Z\omega\,=\,\text{d}(i_Z\omega)+i_Z(\text{d}\omega)\,
=\,\text{d}(\textsc{w}-\text{d}\mathsf{T})+i_Z0\,=\,\text{d}\textsc{w}\ .
\]
\end{proof}

\begin{theorem}\label{thm:conservative-symplectic}
If the mechanical system $(Q,g,F)$ is conservative, its dynamical vector field preserves the symplectic structure of phase space:
\[
L_Z\omega\,=\,0\ .
\]
The converse holds if $Q$ is simply connected.\footnote{It is enough to assume that $H^1(Q,\mathbb R)=0$.}
\end{theorem}
\begin{proof}
If the system is conservative, then $\text{d}\textsc{w}=0$, and Corollary~\ref{cor:lie-work} implies that $L_Z\omega=0$.

If $Q$ is simply connected, so is $TQ$, and hence the $1$-form $\textsc{w}$ is exact if and only if it is closed. By the identity $L_Z\omega=\text{d}\textsc{w}$, this is equivalent to $Z$ preserving the symplectic structure, that is, to $L_Z\omega=0$.
\end{proof}

Every SODE $Z$ on a Riemannian manifold $(Q,g)$ is the dynamical vector field of a mechanical system $(Q,g,F)$, with force intensity defined by $\nabla_Z\mathbb{U}=F$. The theorem says that, when $Q$ is simply connected, this mechanical system is conservative if and only if $Z$ preserves the canonical symplectic structure of $TQ$.

\section{Conservation laws}
\label{sec:conservation-laws}

\subsection*{Restricted symmetries}
Let $(Q,g,F)$ be a mechanical system. Let $\omega$ and $Z$ denote the
symplectic $2$-form and the dynamical vector field on $TQ$, respectively.

\begin{definition}
A hypersurface $H$ (a smooth submanifold of codimension one) of $TQ$, or
of an open subset of $TQ$, is called \emph{special} if the radical of
$\omega|_H$ is $\langle Z\rangle$. In particular, $Z$ is tangent to $H$
and does not vanish on $H$.
\end{definition}

A trajectory $\gamma\colon I\to Q$ of the mechanical system is called
\emph{$H$-restricted} if its velocity lift $\bar\gamma\colon I\to TQ$
takes values in $H$. Since the radical of $\omega|_H$ is
$\langle Z\rangle$, we obtain the following observation.

\begin{observation}\label{obs:restricted-trajectories}
The restriction $\omega|_H$ determines the $H$-restricted trajectories
of the mechanical system, up to reparametrization.
\end{observation}

\begin{definition}
A vector field $D$ tangent to a special hypersurface $H$ is called an
\emph{$H$-restricted infinitesimal symmetry} of the mechanical system if
$L_D(\omega|_H)=0$. It is called an $H$-restricted infinitesimal symmetry
\emph{in the strong sense} if $L_D(\Theta|_H)=0$.
\end{definition}

\begin{proposition}\label{prop:restricted-symmetry-flow}
If $D$ is an $H$-restricted infinitesimal symmetry, then
$L_DZ\in\langle Z\rangle$ on $H$.

Consequently, the flow of $D$ maps $H$-restricted trajectories of the
mechanical system to $H$-restricted trajectories, up to reparametrization
and regarded as curves in phase space $H\subset TQ$.
\end{proposition}
\begin{proof}
Taking the Lie derivative $L_D$ of the equality
$\,0=i_Z(\omega|_H)$ gives
\[
 0\,=\,i_{L_DZ}(\omega|_H)+i_Z\bigl(L_D(\omega|_H)\bigr)\,
   =\,i_{L_DZ}(\omega|_H)+i_Z0\,
   =\,i_{L_DZ}(\omega|_H)\ .
\]
Therefore $L_DZ\in\operatorname{rad}(\omega|_H)=\langle Z\rangle$.
\end{proof}

\begin{definition}
For a vector field $D$ on $H$, the identity
$L_D(\omega|_H)=\text{d}i_D(\omega|_H)$ shows that $D$ is an
$H$-restricted infinitesimal symmetry if and only if the $1$-form
$\iota_D(\omega|_H)$ is closed. If this form is also exact,
$i_D(\omega|_H)=-\text{d}h$ (as happens automatically when $H$ is simply
connected), then the function $h$ on $H$ is called a \emph{Noether
invariant associated with the infinitesimal symmetry $D$}.

If $D$ is an infinitesimal symmetry in the strong sense, we can take
\[
 h:=\,\Theta|_H(D)
\]
without assuming that $H$ is simply connected, because
\[
 \text{d}h\,=\,\text{d}\bigl(\Theta|_H(D)\bigr)\,
   =\,L_D(\Theta|_H)-i_D\text{d}(\Theta|_H)\,
   =\,-i_D(\omega|_H).
\]
This function is given by
\[
 h(v_q)\,=\,g_q\bigl(v_q,\pi_*D_{v_q}\bigr) .
\]
\end{definition}

\begin{theorem}[Noether-Souriau]\label{thm:restricted-noether}
Let $h\in C^\infty(H)$ be the Noether invariant associated with an
$H$-restricted infinitesimal symmetry $D$. Then $h$ is constant along
the $H$-restricted trajectories of the motion: $Z(h)=0$.

Thus, for every $H$-restricted trajectory $\gamma\colon I\to Q$ of the
mechanical system, with velocity $v_t$, the function $h(v_t)$ is
constant along the trajectory.
\end{theorem}
\begin{proof}
\[
 Z(h)\,=\,i_Z\text{d}h\,=\,-i_Zi_D(\omega|_H)\,
      =\,i_Di_Z(\omega|_H)\,=\,i_D0\,=\,0\ .
\]
\end{proof}

\begin{observation}\label{obs:strong-restricted-noether}
If $D$ is an $H$-restricted symmetry in the strong sense, then, for
every $H$-restricted trajectory $\gamma\colon I\to Q$ of the mechanical
system with velocity $v_t$, the function
\[
 h(v_t)\,=\,g_{\gamma(t)}\bigl(v_t,\pi_*D_{v_t}\bigr)
\]
is constant along the trajectory.
\end{observation}

The preceding theorem is an adaptation of a result by Souriau (see \cite{Souriau} Theorems 11.10 and 11.12) to our situation.

\begin{observation}\label{obs:presymplectic-symmetry}
An $H$-restricted infinitesimal symmetry $D$ of the mechanical system is
precisely a locally Hamiltonian vector field for the presymplectic
manifold $(H,\omega|_H)$. If $\iota_D(\omega|_H)=-\text{d}h$ is exact, then
$D$ is a Hamiltonian vector field associated with the function
$h\in C^\infty(H)$.
\end{observation}

\subsection*{Global symmetries}

\begin{definition}
A mechanical system $(Q,g,F)$ is called \emph{quasi-conservative} if
\[
 i_Z\omega\,=\,\textsc{w}-\text{d}\mathsf{T}\,=\,-f\,\text{d}\mathsf{E}
\]
for some smooth functions $f,\mathsf{E}\in C^\infty(TQ)$; the first equality is
the identity established in Theorem~1.1. The function $\mathsf{E}$ is called a
\emph{pseudo-energy}.

If $f=1$, the mechanical system is conservative and $\mathsf{E}$ is its total
energy, up to an additive constant.

Quasi-conservative mechanical systems for which $f$ is nowhere zero
have appeared in the literature under the name \emph{conformally
Hamiltonian systems} (see \cite{Carinena2022,GarciaNaranjo2021,Marle2012}).

Quasi-conservativity implies that, on the open subset of $TQ$ where
$Z$ does not vanish (and hence $f$ does not vanish either), the
distribution defined by the $1$-form
$i_Z\omega=\textsc{w}-\text{d}\mathsf{T}$ is integrable. On this open set,
$\frac1fZ$ is a Hamiltonian vector field.

A quasi-conservative system satisfies the modified Hamilton equation
$$
 i_Z\omega\,=\,-f\,\text{d}\mathsf{E}\ .
$$
\end{definition}

\begin{proposition}\label{prop:pseudoenergy-conservation}
For a quasi-conservative mechanical system, the pseudo-energy is
constant along trajectories of the motion: $Z(\mathsf{E})=0$.
\end{proposition}
\begin{proof}
At a point where $Z$ vanishes, the equality $Z(\mathsf{E})=0$ is immediate. Where
$Z$ does not vanish, neither does the $1$-form
$i_Z\omega=-f\,\text{d}\mathsf{E}$, so $f$ does not vanish. The identity
$0=i_Zi_Z\omega=-fZ(\mathsf{E})$ then implies $Z(\mathsf{E})=0$.
\end{proof}

\begin{observation}\label{obs:pseudoenergy-leaves}
Let $(Q,g,F)$ be a quasi-conservative mechanical system. Write
$T_{\mathrm{reg}}Q$ for the open subset of $TQ$ consisting of regular
points of the dynamical vector field $Z$, that is, points where $Z$
does not vanish. This open subset is dense, since it contains the
complement of the zero section of $TQ\to Q$.
The identity $i_Z\omega=-f\,\text{d}\mathsf{E}$ shows that $\text{d}\mathsf{E}$ does not vanish
where $Z$ does not vanish. Thus, for every $e\in\mathbb R$, the
corresponding nonempty level set
$$
 H_e:=\,\{v_q\in T_{\mathrm{reg}}Q:\mathsf{E}(v_q)=e\}
$$
is a smooth submanifold of codimension one in $T_{\mathrm{reg}}Q$.
Moreover, $H_e$ is a special hypersurface: $Z$ is tangent to it because
$i_Z\text{d}\textsf{E}=Z(\textsf{E})=0$, and
$$
 i_Z(\omega|_{H_e})\,=\,(-f\,\text{d}\mathsf{E})|_{H_e}\,=\,0\ ,
 \qquad \operatorname{rad}(\omega|_{H_e})=\langle Z\rangle.
$$
In conclusion, the pseudo-energy level hypersurfaces define a
foliation of $T_{\mathrm{reg}}Q$ by special hypersurfaces.

In this setting, Observation~\ref{obs:restricted-trajectories} can be
restated as follows: the restriction $\omega|_{H_e}$ determines the
trajectories of pseudo-energy $\mathsf{E}=e$, up to reparametrization.
\end{observation}

\begin{definition}
For a quasi-conservative mechanical system $(Q,g,F)$ with
pseudo-energy $\mathsf{E}$, a \emph{global infinitesimal symmetry} is a vector
field $D$ on $TQ$ satisfying $L_D\omega=0$ and $D(\mathsf{E})=0$. It is a global
infinitesimal symmetry \emph{in the strong sense} if
$L_D\Theta=0$ and $D(\mathsf{E})=0$.

Occasionally, infinitesimal symmetries defined only on an open subset
of $TQ$ will be considered.
\end{definition}

Let $D$ be a global infinitesimal symmetry. The condition $D(\mathsf{E})=0$
implies that $D$ is tangent to every hypersurface $H_e$. Thus every
global symmetry is also an $H_e$-restricted symmetry for every
$e\in\mathbb R$. The converse statement is false.

The definition of a Noether invariant extends to global symmetries.

\begin{definition}
If $D$ is a global infinitesimal symmetry, the $1$-form
$i_D\omega$ is closed, since
$$
 \text{d}i_D\omega\,=\,L_D\omega-i_D\text{d}\omega\,=\,0-i_D0\,=\,0\ .
$$
If it is also exact, $i_D\omega=-\text{d}h$, the function $h$ on $TQ$ is
called a \emph{Noether invariant associated with the symmetry $D$}.

If $D$ is a global infinitesimal symmetry in the strong sense, we can
take $h:=\Theta(D)$.
\end{definition}

Global symmetries preserve the line field of the motion at regular
points and satisfy Noether's theorem. Their proofs are similar to
those for $H$-restricted symmetries.

\begin{proposition}\label{prop:global-symmetry-flow}
Let $D$ be a global infinitesimal symmetry of a quasi-conservative
mechanical system. Then $L_DZ\in\langle Z\rangle$ on
$T_{\mathrm{reg}}Q$.

Moreover, if the system is conservative, then $L_DZ=0$ on $TQ$; that is, the
flow of $D$ maps trajectories of the motion in phase space to
trajectories of the motion while preserving their parametrization.
\end{proposition}
\begin{proof}
Taking the Lie derivative $L_D$ of the equality
$i_Z\omega=-f\,\text{d}\mathsf{E}$ gives
$$
 i_{(L_DZ)}\omega+\underbrace{i_Z(L_D\omega)}_{0}\,
 =\,-D(f)\,\text{d}\mathsf{E}-f\,\underbrace{\text{d}\bigl(D(\mathsf{E})\bigr)}_{0}\ .
$$
Thus $i_{(L_D)Z}\omega=-D(f)\,\text{d}\mathsf{E}$. Comparing this with
$i_Z\omega=-f\,\text{d}\mathsf{E}$ yields
$$
 L_DZ=\frac{D(f)}{f}Z
$$
wherever $Z$ does not vanish, and hence $f$ does not vanish.

For a conservative system, choosing the total energy as $\mathsf{E}$ gives
$f=1$, and therefore $L_DZ=0$.
\end{proof}

\begin{theorem}[Noether]\label{thm:global-noether}
Let $(Q,g,F)$ be a quasi-conservative mechanical system. For every
$h\in C^\infty(TQ)$, the following conditions are equivalent:
\begin{enumerate}
 \item[(a)] $h$ is a Noether invariant: the vector field $D$ defined by
 $\,i_D\omega=-\mathrm{d}h\,$ is a global infinitesimal symmetry.
 \item[(b)] $h$ is constant along trajectories of the motion:
 $Z(h)=0$.
 \item[(c)] $\{h,\mathsf{E}\}=0$.
\end{enumerate}
\end{theorem}
\begin{proof}
For $(a)\Rightarrow(b)$,
$$
 \begin{aligned}
 Z(h)&\,=\,i_Z\text{d}h\,=\,-i_Zi_D\omega\,
       =\,i_Di_Z\omega\\
     &\,=\,-i_D(f\,\text{d}\mathsf{E})\,=\,-fD(\mathsf{E})\,=\,0\ .
 \end{aligned}
$$
For $(b)\Rightarrow(a)$, we verify that $D$ is an infinitesimal symmetry.
First,
$$
 L_D\omega\,=\,\text{d}i_D\omega+i_D\text{d}\omega\,
           =\,-\text{d}(\text{d}h)+i_D0\,=\,0\ .
$$
Second, on $T_{\mathrm{reg}}Q$,
$$
 \begin{aligned}
 D(\mathsf{E})&\,=\,i_D\text{d}\mathsf{E}\,
      =\,i_D\left(-\frac1fi_Z\omega\right)\,
      =\,\frac1fi_Zi_D\omega\\
     &\,=\,-\frac1fi_Z\text{d}h\,=\,-\frac1fZ(h)\,=\,0\ .
 \end{aligned}
$$
By continuity, $D(\mathsf{E})=0$ on all of $TQ$.

For $(a)\Rightarrow(c)$, we have $\{h,\mathsf{E}\}:=D(\mathsf{E})=0$.
Conversely, $(c)$ gives $0=\{h,\mathsf{E}\}=D(\mathsf{E})$, while
$L_D\omega=\text{d}i_D\omega+i_D\text{d}\omega=-\text{d}(\text{d}h)+i_D0=0$. Hence $D$ is
an infinitesimal symmetry, proving $(c)\Rightarrow(a)$.

\end{proof}

\begin{observation}\label{obs:symmetry-integral-correspondence}
If $Q$ is connected and simply connected, then so is $TQ$, and every
global infinitesimal symmetry produces a Noether invariant. There is
therefore a bijective correspondence
$$
\left\{\begin{aligned}
\text{global infinitesimal symmetries}\\
\text{of the mechanical system}\qquad
\end{aligned}\right\} \longleftrightarrow 
\left\{\begin{aligned}
\text{first integrals of the}\quad\quad\\
\text{dynamical vector field }Z
\end{aligned}\right\}/\mathbb{R}\ ,\quad D \longleftrightarrow [h]
$$
defined by $i_D\omega=-\text{d}h$.
\end{observation}

\subsection*{Killing vector fields}
Every vector field $D$ on $Q$ has a natural lift $\widetilde D$ to
$TQ$, also called its \emph{complete lift}. If $\{\phi_t\}$ is the
local flow of $D$, this lift is defined as the infinitesimal generator
of the tangent maps $\{{\phi_t}_*\}$. In local coordinates
$(q_1,\ldots,q_n,v_1,\ldots,v_n)$ on $TQ$, if
\[
 D=\sum_i h_i\partial_{q_i},
\]
then
\begin{equation}\label{eq:complete-lift}
 \widetilde D
 =\sum_i h_i\partial_{q_i}
  +\sum_{i,j}\frac{\partial h_i}{\partial q_j}\,v_j\partial_{v_i}.
 \tag{2}
\end{equation}
Notice that $\pi_*(\widetilde D)=D$.

\begin{proposition}\label{prop:killing-symmetry}
Let $(Q,g,F)$ be a conservative mechanical system. Let $K$ be a
Killing vector field on $(Q,g)$ that preserves the potential:
$K(\mathsf{U})=0$. Then its natural lift $\widetilde K$ to $TQ$ is an
infinitesimal symmetry of the mechanical system in the strong sense.
\end{proposition}
\begin{proof}
If $\sigma\colon Q\to Q$ is a diffeomorphism, its tangent map
$\sigma_*\colon TQ\to TQ$ is also a diffeomorphism. If $\sigma$ is an
isometry, then $\sigma_*$ preserves the Liouville form and the kinetic
energy:
\[
 (\sigma_*)^*\Theta\,=\,\Theta\ ,
 \qquad \mathsf{T}\circ \sigma_*\,=\,\mathsf{T}.
\]
If $\{\phi_t\}$ is the local one-parameter group of a vector field
$K$ on $Q$, then $\{\phi_{t\,*}\}$ is the local one-parameter group of
$\widetilde K$. When $K$ is Killing, the maps $\phi_t$ are isometries.
Consequently, the maps $\phi_{t\,*}$ preserve the Liouville form and the
kinetic energy, giving
\[
 L_{\widetilde K}\Theta\,=\,0\ ,
 \qquad \widetilde K(\mathsf{T})\,=\,0\ .
\]
Hence
\[
 \widetilde K(E)\,=\,\widetilde K(\mathsf{T})
                +\widetilde K(\mathsf{U})\,=\,0\ ,
\]
where functions on $Q$ are understood to be pulled back by $\pi$.
Thus $\widetilde K$ is an infinitesimal symmetry in the strong sense.

\end{proof}

\begin{corollary}\label{cor:killing-integral}
Let $(Q,g,F)$ be a conservative mechanical system. If $K$ is a
Killing vector field on $(Q,g)$ preserving the potential energy,
$K(\mathsf{U})=0$, then, along every solution curve
$\gamma\colon I\to Q$ of the mechanical system with velocity $v$, the
scalar product $g(v,K)$ is constant.
\end{corollary}
\begin{proof}
This follows from Proposition~\ref{prop:killing-symmetry} and
Observation~\ref{obs:strong-restricted-noether}.

\end{proof}

The preceding two results extend to quasi-conservative systems, with
a similar proof.

\begin{proposition}\label{prop:quasiconservative-killing}
Let $(Q,g,F)$ be a quasi-conservative mechanical system. Let $K$ be a
Killing vector field on $(Q,g)$ whose natural lift preserves the
pseudo-energy: $\widetilde K(E)=0$. Then $\widetilde K$ is an
infinitesimal symmetry of the mechanical system in the strong sense.

Consequently, along every solution curve $\gamma\colon I\to Q$ with
velocity $v$, the scalar product $g(v,K)$ is constant.
\end{proposition}

\section{Jacobi's theorem}

The classical Jacobi theorem states that the trajectories of a conservative
mechanical system with prescribed total energy coincide, up to
reparametrization, with the geodesics of a Riemannian metric. We shall give a
slightly more general version of this theorem.

\begin{theorem}\label{thm:jacobi-general}
Let $(Q,g,F)$ be a mechanical system, and let $\alpha\in C^\infty(Q)$ be a
positive function. If
\[
 H:=\{v_q\in TQ:\mathsf{T}(v_q)=\alpha(q)\}
\]
is a special hypersurface, then the $H$-restricted trajectories of the
mechanical system coincide, up to reparametrization, with the nonconstant geodesics of
the Riemannian metric $\bar g=\alpha g$.
\end{theorem}

\begin{proof}
The geodesics of the metric $\bar g:=\alpha g$ are the solutions of the
mechanical system $(Q,\bar g,\bar F=0)$. Its Liouville form is
$\bar\Theta=\alpha\Theta$, and its symplectic form is
$\bar\omega=\text{d}\bar\Theta$. The total energy of this mechanical system is
$$
 \bar {\mathsf{E}}(v_q)\,=\,\bar{\mathsf{T}}(v_q)\,=\,\alpha(q)\mathsf{T}(v_q)\ .
$$
Thus its energy hypersurface $\bar{\mathsf{E}}=1$ is
$$
 \bar H\,=\,\{v_q\in TQ:\alpha(q)\mathsf{T}(v_q)=1\}\ .
$$
It is immediate that the diffeomorphism
$$
 \varphi:TQ\longrightarrow TQ,\qquad v_q\longmapsto\alpha(q)v_q
$$
maps $\bar H$ onto $H$:
$$
 \varphi:\bar H\longrightarrow H,\qquad v_q\longmapsto\alpha(q)v_q.
$$
Moreover, $\varphi^*\Theta=\alpha\Theta=\bar\Theta$. Differentiation
therefore gives
\[
 \varphi^*(\omega|_H)\,=\,\bar\omega|_{\bar H}\ .
\]
By Observation~2.1, this implies that the diffeomorphism
$\varphi:\bar H\to H$ satisfies
$\varphi_*\langle\bar Z\rangle=\langle Z\rangle$. Consequently,
$\varphi$ maps the trajectories of $\bar Z|_{\bar H}$, which are lifts of
geodesics of $\bar g$, onto the trajectories of $Z|_H$, which are lifts of
$H$-restricted solutions of the original mechanical system, up to
reparametrization. Projecting both kinds of curves onto the base $Q$, and
using the commutative triangle
$$
 \begin{tikzcd}
  \bar H \arrow[rr,"\varphi"] \arrow[dr,"\pi"'] &&
  H \arrow[dl,"\pi"]\\
  &Q&
 \end{tikzcd}
$$
shows that the geodesics of $\bar g$ coincide with the $H$-restricted
solutions of $(Q,g,F)$.

\end{proof}

\begin{corollary}\label{cor:jacobi-killing}
With the notation of Theorem~3.1, every Killing field $K$ of $\bar g$
induces a strong $H$-restricted infinitesimal symmetry $\mathbb K$, whose
Noether invariant is
$$
 h_{\mathbb K}(v_q)\,=\,g(v_q,K_q),\qquad v_q\in H.
$$

The mapping $K\mapsto\mathbb{K}$ establishes an isomorphism between the Lie algebra of Killing vectors of $\bar g$ and an algebra of $H$-restricted symmetries of $(Q,g,F)$.

\end{corollary}

\begin{proof}
Let $\widetilde K$ be the lift of $K$ to $TQ$. This vector field is tangent
to $\bar H$ and is a strong $\bar H$-restricted infinitesimal symmetry of
the mechanical system $(Q,\bar g,\bar F=0)$. Since the diffeomorphism
$\varphi:\bar H\to H$ satisfies
$\varphi^*(\Theta|_H)=\bar\Theta|_{\bar H}$, it follows that
\[
 \mathbb K:=\varphi_*(\widetilde K|_{\bar H})
\]
is a strong restricted infinitesimal symmetry on $H$. The corresponding
Noether invariants satisfy
$h_{\mathbb K}=h_{\widetilde K}\circ\varphi^{-1}$. Hence
\begin{align*}
 h_{\mathbb K}(v_q)
 &\,=\,h_{\widetilde K}(\varphi^{-1}(v_q))
  \,=\,h_{\widetilde K}\!\left(\frac{v_q}{\alpha(q)}\right)\\
 &\,=\,\bar g\!\left(\frac{v_q}{\alpha(q)},K_q\right)
  \,=\,g(v_q,K_q).
\end{align*}
\end{proof}

The classical Jacobi theorem is a particular case of the preceding theorem.

\begin{theorem}[Jacobi]\label{thm:jacobi-classical}
Let $(Q,g,F=-\operatorname{grad}\mathsf{U})$ be a conservative mechanical
system. Given $e\in\mathbb R$, let $Q_0$ be the open subset of $Q$ on which
$e-\mathsf{U}>0$. The solutions in $Q_0$ with total energy $e$ coincide,
up to reparametrization, with the nonconstant geodesics of the Riemannian metric
\[
 \bar g\,=\,(e-\mathsf{U})g|_{Q_0}\ .
\]
Moreover, every Killing vector field $K$ of the Jacobi metric
$\bar g=(e-\mathsf{U})g$ induces an $H_e$-restricted infinitesimal symmetry
$\mathbb K$ of the conservative mechanical system $(Q,g,\mathsf{U})$,
whose Noether invariant is
$$
 h_{\mathbb K}(v_q)\,=\,g(v_q,K_q)\ ,\qquad v_q\in H_e\ .
$$
\end{theorem}

\begin{proof}
The hypersurface $H_e$ in $TQ_0$ defined by $\mathsf{E}=e$ is
\begin{align*}
 H_e
 &\,=\,\{v_q\in TQ_0:E(v_q)=\mathsf{T}(v_q)+\mathsf{U}(q)=e\}\\
 &\,=\,\{v_q\in TQ_0:\mathsf{T}(v_q)=\alpha(q)\},
\end{align*}
where $\alpha(q):=e-\mathsf{U}(q)$.
For every $v_q\in H_e$, one has
$\mathsf{T}(v_q)=e-\mathsf{U}(q)>0$. Thus $v_q\ne0$ and consequently
$Z_{v_q}\ne0$. Hamilton's equation implies that $\text{d}\mathsf{E}$ does not vanish at
any point of $H_e$. Therefore $H_e$ is a smooth hypersurface; moreover,
it is special, as observed in \ref{obs:pseudoenergy-leaves}. Apply Theorem~3.1 and Corollary~3.2.

\end{proof}

The generalized Jacobi theorem \ref{thm:jacobi-general} can be derived from the classical version by applying the following result.

\begin{proposition} Let $(Q,g,F)$ be a mechanical system, and let $\alpha\in C^\infty(Q)$ be a
positive function. The hypersurface
$$
H\,=\,\{v_q\in TQ:\mathsf{T}(v_q)=\alpha(q)\}
$$
is special if and only if
$$
F_{v_q}\,=\,(\mathrm{grad}\,\alpha)_q
\qquad\text{for every }v_q\in H.
$$
Thus, on $H$, the system coincides with the conservative mechanical
system with potential $-\alpha$, restricted to its zero-energy
hypersurface.
\end{proposition}
\begin{proof}  Let $V$ be the vertical Euler vector field on $TQ$.
Since
$$V(\mathsf{T}-\alpha)\,=\,2\mathsf{T}-0\,\overset{\text{on }H}{=}\,2\alpha\,>\,0\ ,
$$
we obtain that $\text{d}(\mathsf{T}-\alpha)$ does not vanish at any point of $H$; hence, the hypersurface $H=(\mathsf{T}-\alpha)^{-1}(0)$  is regular.

$(\Rightarrow)$. If $H$ is special, then $(i_Z\omega)_{|H}=0$, that is to say, $i_Z\omega$ is proportional to $\text{d}(\mathsf{T}-\alpha)$ at the points of $H$, so that
$$
(\textsc{w}-\text{d}\mathsf{T})_{v_q}\,
=\,(i_Z\omega)_{v_q}
\,=\, \lambda(v_q)\,\text{d}(\mathsf{T}-\alpha)_{v_q},
\qquad v_q\in H.
$$
for some $\lambda\in\mathcal{C}^\infty(H)$. Evaluating it on $V$, and using $\textsc{w}(V)=0$, gives
$$
-2\mathsf{T}\,=\,\lambda\,2\mathsf{T},
$$
hence $\lambda=-1$. Therefore
$$
\textsc{w}_{v_q}\,=\,(\text{d}\alpha)_{v_q} ,
\qquad v_q\in H.
$$
By the definition of the work form, this is equivalent to
$F_{v_q}=(\text{grad}\,\alpha)_q$.

($\Leftarrow)$. The identity $\textsc{w}_{v_q}\,=\,(\text{d}\alpha)_{v_q}$ for every $v_q\in H$ implies
$$
(i_Z\omega)_{v_q}\,=\,(\textsc{w}-\text{d}\mathsf{T})_{v_q}\,
=\,-\text{d}(\mathsf{T}-\alpha)_{v_q}\qquad\text{for every }v_q\in H.
$$
Antisymmetry gives $i_Z\text{d}(\mathsf{T}-\alpha)=0$ on $H$, so $Z$ is tangent
to $H$ and spans the kernel of $\omega|_H$.
Indeed, this kernel is one-dimensional, and $Z$ does not vanish
on $H$ because $\pi_*Z_{v_q}=v_q\ne0$.
Thus $H$ is special.

\end{proof}


Under certain general hypotheses, the special hypersurfaces of a
mechanical system, when they exist, have the form used in Theorem~3.1.

\begin{proposition}\label{prop:spherical-special-hypersurfaces}
Let $H\subset TQ$ be a special hypersurface of a mechanical system
$(Q,g,F)$, with $\mathrm{dim}\, Q\geq 2$. Suppose that $\pi:H\to Q$ is a surjective submersion with
connected, nonempty fibers that are closed in $TQ$. Then
\[
 H\,=\,\{v_q\in TQ:\mathsf{T}(v_q)=\alpha(q)\}
\]
for some positive function $\alpha\in C^\infty(Q)$.
\end{proposition}

\begin{proof}
Let $D$ be a vector field on $H$ satisfying $\pi_*D=0$, that is, tangent
to the fibers of $\pi:H\to Q$. Notice that
$\textsc{w}(D)=g(F,\pi_*D)=0$. Since $i_Z(\omega|_H)=0$, on $H$ we have
\[
 0\,=\,i_D(i_Z\omega)\,=\,i_D(\textsc{w}-\text{d}\mathsf{T})\,
  =\,\textsc{w}(D)-D\mathsf{T}\,=\,-D\mathsf{T}\ .
\]
Thus $\mathsf{T}$ is constant on the fibers of $\pi:H\to Q$:
$\mathsf{T}(v_q)=\alpha(q)$ for every $v_q\in H$, for some
$\alpha\in C^\infty(Q)$. Consequently, the fiber over $q\in Q$ satisfies
\[
 (\pi|_H)^{-1}(q)
 \,=\,\{v_q\in T_qQ:v_q\in H\}\,
 \subseteq\,\{v_q\in T_qQ:\mathsf{T}(v_q)=\alpha(q)\}.
\]
The set on the right is the sphere of radius $\sqrt{2\alpha(q)}$ in
$T_qQ$, since $\mathsf{T}(v_q)=\tfrac12|v_q|^2$. The inclusion is an
equality because $(\pi|_H)^{-1}(q)$ is a closed submanifold of codimension
one in $T_qQ$.

\end{proof}

\section{Central forces and hodographs}

Consider Euclidean space $(\mathbb R^n,g=\sum_i \text{d}x_i^2)$, and let
$(x_1,\ldots,x_n,v_1,\ldots,v_n)$ be the corresponding coordinates on
$T\mathbb R^n$. These coordinates define a diffeomorphism
\[
 T\mathbb R^n\longrightarrow\mathbb R^n\times\mathbb R^n,
 \qquad v_x\longmapsto(x,v_x).
\]

\begin{definition}
Let $J$ denote the following diffeomorphism:
\[
 J:T\mathbb R^n=\mathbb R^n\times\mathbb R^n
 \longrightarrow\mathbb R^n\times\mathbb R^n=T\mathbb R^n,
 \qquad J(x,v)=(v,x).
\]
Notice that $J$ preserves the symplectic structure
$\omega=\sum_i \text{d}v_i\wedge \text{d}x_i$ up to a sign:\footnote{A slight
modification of the definition, $J(x,v)=(v,-x)$, would make $J$ a
symplectomorphism.}
\[
 J^*\omega\,=\,-\omega.
\]
\end{definition}

\begin{definition}
The \emph{hodograph} of a curve $x(t)$ in $\mathbb R^n$ is the curve
$x'(t)$.
\end{definition}

Let $x(t)$ be a curve in $\mathbb R^n$, and consider its lift
$(x(t),x'(t))$ to $T\mathbb R^n$. Its image under $J$ is
$(x'(t),x(t))$. This image curve is itself the lift of a curve in
$\mathbb R^n$---the hodograph, possibly with a different
parametrization---if and only if $x''(t)$ is proportional to $x(t)$ with
a nonzero proportionality factor. This occurs, for instance, when
$x(t)$ is the trajectory of a particle moving under a central force.
Consider the following definition.

\begin{definition}
Let $h$ be a smooth, nowhere-vanishing function on $T\mathbb R^n$, or on
an open subset of $T\mathbb R^n$. Let $Z$ and $Z^*$ be the respective
motion vector fields of the second-order differential equations
\[
 x''=h(x,x')x\ ,
 \qquad
 x''\,=\,\frac{1}{h(x',x)}x\ .
\]
We say that these second-order differential equations are \emph{mutually
dual}.

The vector fields $Z$ and $Z^*$ are the motion vector fields of two
mechanical systems on $(\mathbb R^n,g=\sum_i \text{d}x_i^2)$, with respective
central forces
\[
 F_{(x,v)}\,=\,h(x,v)x\ ,
 \qquad
 F^*_{(x,v)}\,=\,\frac{1}{h(v,x)}x\ .
\]
We also call these mechanical systems \emph{mutually dual}.
\end{definition}

Notice that multiplying the Euclidean metric by a positive constant
changes neither $J$ nor the duality of the second-order differential
equations and the mechanical systems.

\begin{theorem}\label{thm:central-duality}
With the preceding notation:
\begin{enumerate}
 \item[(a)] The diffeomorphism $J:T\mathbb R^n\to T\mathbb R^n$ maps
 the second-order differential equation $Z$ onto its dual $Z^*$ up to
 a scalar factor:
 \[
  J_*\langle Z\rangle\,=\,\langle Z^*\rangle\ .
 \]
 \item[(b)] The solutions of $Z^*$ are the reparametrized hodographs of
 the solutions of $Z$.
 \item[(c)] The diffeomorphism $J$ maps special hypersurfaces of $Z$
 onto special hypersurfaces of the dual equation $Z^*$.
 \item[(d)] If the original mechanical system is quasi-conservative
 with pseudo-energy $\mathsf{E}$, then the dual system is also
 quasi-conservative, with pseudo-energy
 $$
  \mathsf{E}^*:=\,\mathsf{E}\circ J^{-1}\ .
$$
\end{enumerate}
\end{theorem}

\begin{proof}
(a) We have
\[
 Z\,=\,\sum_i v_i\partial_{x_i}+h(x,v)\sum_i x_i\partial_{v_i}.
\]
Applying $J_*$ immediately gives
\begin{align*}
 J_*Z
 &\,=\,\sum_i x_i\partial_{v_i}+h(v,x)\sum_i v_i\partial_{x_i}\\
 &\,=\,h(v,x)\left(\sum_i v_i\partial_{x_i}
    +\frac{1}{h(v,x)}\sum_i x_i\partial_{v_i}\right)
  \,=\, h(v,x)Z^*.
\end{align*}

(b) Let $x(t)$ be a solution of $Z$, so that its lift
$(x(t),x'(t))$ to $T\mathbb R^n$ is an integral curve of $Z$. The image
curve
\[
 J(x(t),x'(t))=(x'(t),x(t))
\]
is an integral curve of $J_*Z\,=\,h(v,x)Z^*$. After reparametrization, it is
an integral curve of $Z^*$ and therefore the lift of a curve $y(s)$:
\[
 (x'(t),x(t))\,=\,(y(s),y'(s))\ .
\]
Consequently, the solutions $y(s)$ of $Z^*$ are the reparametrized
hodographs $y(s)=x'(t)$ of the solutions $x(t)$ of $Z$.

(c) This follows from $J^*\omega=-\omega$ and
$J_*\langle Z\rangle=\langle Z^*\rangle$.

(d) Write $J_*Z=\bar f Z^*$, where $\bar f=h\circ J^{-1}$ is nowhere
zero. Starting from $i_Z\omega=-f\,\text{d}\mathsf{E}$, pushing forward by $J$ gives
$$
 i_{\bar f Z^*}(-\omega)
 \,=\,-(f\circ J^{-1})\,\text{d}(\mathsf{E}\circ J^{-1}).
$$
Thus
$$
 i_{Z^*}\omega
 \,=\,\frac{f\circ J^{-1}}{\bar f}\,\text{d}\mathsf{E}^*
 \,=\,-f^*\,\text{d}\mathsf{E}^*\ ,
 \qquad
 f^*:=\,-\frac{f\circ J^{-1}}{\bar f}\ ,
 \qquad \mathsf{E}^*:=\,\mathsf{E}\circ J^{-1}\ .
$$
\end{proof}

\refstepcounter{thm}\label{remark:dual-noether}
\noindent\textbf{\thethm\  Remark.}
Let $H_e$ be a level hypersurface $\mathsf{E}=e$ of the pseudo-energy of a
quasi-conservative system. Then $H_e^*:=J(H_e)$ is a level hypersurface
of the pseudo-energy $\mathsf{E}^*:=\mathsf{E}\circ J^{-1}$ of the dual system. If $D$ is
an $H_e$-restricted symmetry of the original system, then $J_*D$ is an
$H_e^*$-restricted symmetry of the dual system, and the corresponding
Noether invariants $h$ and $h^*$ satisfy
$$
 h^*\,=\,-h\circ J^{-1}\ ,
$$
where the sign is a consequence of the anti-symplectic character of $J$.\medskip

Let us consider the case of a conservative mechanical system
with a central force.

\begin{theorem}\label{4.3} Let $(\mathbb R^n\setminus\{0\},g=\sum_i \mathrm{d}x_i^2,F)$ be a conservative
mechanical system with central force and total energy
$$F_x\,=\,-\operatorname{grad}\mathsf{U}(|x|)\quad ,\qquad \mathsf{E}(x,v)\,=\,\frac{1}{2}|v|^2+\mathcal{U}(|x|)\ ,$$ where
$\mathsf{U}'(r)\ne0$ for every $r>0$.\smallskip

(a) The dual mechanical system has central force and pseudo-energy $$F^*_{(x,v)}\,=\,-\frac{|v|}{\mathsf{U}'(|v|)}x\quad ,\qquad   \mathsf{E}^*(x,v)\,=\,\tfrac12|x|^2+\mathsf{U}(|v|)\ .$$

(b) The hodographs of the solutions of the original system with total energy $e$ are the reparametrized solutions  of the dual system with pseudo-energy $e$.\smallskip

(c) The solutions of the dual mechanical system with pseudo-energy $e$ are the nonconstant reparametrized geodesics of the Riemannian metric 
$$
 \bar g_e\,=\,\mathsf V^2\!\left(e-\tfrac12|x|^2\right)g\ ,$$
where $\mathsf V$ is the inverse function of $\,\mathsf{U}$.\smallskip

(d) Every  Killing vector field $K^e$ of $\bar g_e$ induces an $H^*_e$-restricted symmetry $\mathbb{K}^{e\,*}$ of the dual problem and an $H_e$-restricted symmetry $\mathbb{K}^e$ of the original mechanical system.
\end{theorem}

\refstepcounter{thm}\label{note}
\noindent\textbf{\thethm\  Note.}
The condition $\mathsf{U}'(r)\ne0$ for all $r>0$ implies that
$\mathsf{U}$ is a diffeomorphism from $(0,\infty)$ onto an interval
$(a,b)$, where
$a=\inf_{r>0}\mathsf{U}(r)$ and
$b=\sup_{r>0}\mathsf{U}(r)$. Its inverse is the diffeomorphism
$\mathsf V:(a,b)\to(0,\infty)$.
The metric $\bar g_e$ is defined on the open subset of $\mathbb R^n$
where $\mathsf V(e-\tfrac12|x|^2)$ is defined, namely
$$
 2e-2b\,<\,|x|^2\,<\,2e-2a\ .
$$

\begin{proof} $(a)$. The trajectories of the original system solve the second-order equation
$$
 x''\,=\,-\frac{\mathsf{U}'(|x|)}{|x|}x\ ,$$
hence the dual second-order equation is
$$
 x''\,=\,-\frac{|x'|}{\mathsf{U}'(|x'|)}x\ .
$$
Thus the dual mechanical system has central force
$$
 F^*_{(x,v)}\,=\,-\frac{|v|}{\mathsf{U}'(|v|)}x\ ,
$$
 with pseudo-energy
$$
 \mathsf{E}^*\,=\,\mathsf{E}\circ J^{-1}\ ,\qquad
 \mathsf{E}^*(x,v)\,=\,\tfrac12|x|^2+\mathsf{U}(|v|)\ .
$$

$(b)$. By Theorem~4.1(b) and Remark~4.2, the required hodographs are the
reparametrized solutions of the dual mechanical system with
pseudo-energy $\mathsf{E}^*=e$.\smallskip

$(c)$.  For the pseudo-energy
$$
 \mathsf{E}^*\,=\,\tfrac12|x|^2+\mathsf{U}(|v|)$$
 the level hypersurface $H_e^*$ has equation
$$
 \tfrac12|x|^2+\mathsf{U}(|v|)\,=\,e\ ,
$$
or equivalently
$$
 |v|\,=\,\mathsf V\!\left(e-\tfrac12|x|^2\right)\ .
$$
Consequently,
$$
 H_e^*\,=\,\left\{v_x\in T\mathbb R^n:
 \mathsf{T}(v_x)=\tfrac12\mathsf V^2
       \!\left(e-\tfrac12|x|^2\right)\right\}\ ,
$$
with $x$ in the domain of the displayed expression.
Apply the generalized Jacobi theorem, Theorem~3.1. The factor
$\tfrac12$ has been omitted from the metric because it does not affect
the geodesics.\smallskip

(d) By Corollary \ref{cor:jacobi-killing}, every Killing vector $K^e$ induces an $H_e^*$-restricted symmetry $\mathbb{K}^{e\,*}$ of the dual problem. Since $J\colon H_e\to H^*_e$ is an anti-symplectic diffeomorphism, we have that $\mathbb{K}^e:=J^{-1}_*\mathbb{K}^{e\,*}$ is an $H_e$-restricted symmetry.

\end{proof}

\section{Application to the Kepler problem}
Consider the motion of a point particle of mass $m>0$ in Euclidean space
$(\mathbb R^n,\sum_i \text{d}x_i^2)$, ($n\geq 2$), under the central force intensity
$F_x=-kr^{-3}x$, where $k$ is a positive constant and $r=|x|$ is the
Euclidean distance from the origin. The corresponding SODE is
\[
 x''\,=\,-k\frac{x}{|x|^3}\ .
\]
The dual SODE and the dual force intensity are
\[
 x''\,=\,-\frac1k|x'|^3x\ ,
 \qquad F^*_{(x,v)}\,=\,-\frac1k|v|^3x\ .
\]
The configuration space of the Kepler problem is the manifold
$\mathbb R^n\setminus\{0\}$, endowed with the kinetic metric
$g=m\sum_i \text{d}x_i^2$. The potential energy is $\mathsf{U}=-km/|x|$, and the
total energy is
$$
 \mathsf{E}\,=\,\tfrac12m|v|^2-\frac{km}{|x|}\ .
$$
The dual problem has pseudo-energy
$\mathsf{E}^*=\tfrac12m|x|^2-km/|v|$.

The kinetic metric is ``Euclidean'': $g=\sum_i \text{d}\bar x_i^2$, where
$\bar x_i=\sqrt m\,x_i$. The distance from the origin with respect to this
metric is $\bar r=\sqrt m\,r$. In terms of this distance the potential
energy is $\mathsf{U}=-km/r=-km^{3/2}/\bar r$. Its inverse has the same
functional form: $\mathsf{V}(t)=-km^{3/2}/t$ for $t<0$.

Let us compute the metric $\bar g_e$ of Theorem~4.3 in the case of the Kepler problem,
\begin{align*}
 \bar g_e
 &\,=\,\mathsf V^2\!\left(e-\tfrac12\bar r^2\right)g
 \,=\,k^2m^3\frac{1}{(e-\tfrac12\bar r^2)^2}\,g\\
 &\,=\,k^2m^4\frac{1}{(e-\tfrac12mr^2)^2}\sum_i \text{d}x_i^2
 \,=\,k^2m^2\frac4{(r^2-2e/m)^2}\sum_i \text{d}x_i^2\ .
\end{align*}
By Note \ref{note}  this metric is defined
on the open subset of $\mathbb R^n$ where $2e<\bar r^2$, or equivalently,
$r^2-2e/m>0$.

Since two Riemannian metrics that differ by a constant positive factor
have the same geodesics, we may omit the initial factor $k^2m^2$ and take
$$
 \bar g_e\,=\,\frac4{(r^2-2e/m)^2}\sum_i \text{d}x_i^2
 \,=\,\frac4{(r^2+\mathsf{K})^2}\sum_i \text{d}x_i^2\ ,
 \qquad \mathsf{K}=-\frac{2e}{m}\ .
$$
This metric has constant curvature $\mathsf{K}$ (see Observation~5.2 below).
Applying Theorem \ref{4.3} we obtain the following result.

\begin{theorem}\label{thm:kepler-hodographs}
Let $(\mathbb R^n\setminus\{0\},g=m\sum_i \mathrm{d}x_i^2,F_x=-kr^{-3}x)$
be the mechanical system of the Kepler problem. 

(a) The hodographs of its
trajectories of energy $e$ are the reparametrized nonconstant geodesics of the
Riemannian metric of constant curvature $\mathsf{K}=-2e/m$,
$$
 \bar g_e\,=\,\frac4{(|x|^2-2e/m)^2}\sum_i \mathrm{d}x_i^2\,
 =\,\frac4{(|x|^2+\mathsf{K})^2}\sum_i \mathrm{d}x_i^2,
$$
on the region $|x|^2-2e/m=|x|^2+\mathsf{K}>0$.

(b) Every Killing vector $K^e$ of $\,\bar g_e$ defines an $H_e$-restricted symmetry $\mathbb{K}^e$ of the Kepler problem.
\end{theorem}

\begin{observation}\label{obs:curvature}
The sectional curvature $\mathsf{K}_{ij}$ of the plane
$\langle\partial_{x_i},\partial_{x_j}\rangle$ for the Riemannian metric
$$
 g\,=\,\frac{\sum_i \text{d}x_i^2}{q(x_1,\ldots,x_n)^2}\ ,\qquad (n\geq 2)
$$
is
$$
\mathsf{K}_{ij}\,=\,(q_{ii}+q_{jj})q-|\operatorname{grad}q|^2\ ,
 \qquad q_{ii}\,=\,\frac{\partial^2q}{\partial x_i^2}\ ,
$$
where the gradient and norm on the right are Euclidean. In particular,
the metric
$$
 g\,=\,\frac{4\sum_i \text{d}x_i^2}{(a\sum_i x_i^2+b)^2}\ ,\qquad a,b\in\mathbb R\ ,
$$
has constant curvature $\mathsf{K}=ab$ wherever its denominator does not vanish.
It follows that the metric $\bar g_e$ in Theorem~5.1 has curvature
$\mathsf{K}=-2e/m$.
\end{observation}

\subsection*{Geometry of the metrics in Theorem 5.1}
We shall see that the metrics $\bar g_e$ in Theorem~5.1 are almost
complete: adding a single point at infinity produces a complete
geometry. According to the value of $e$, we describe their relationship
with the complete space-form models $\mathbb{S}^n$, $\mathbb R^n$, and $\mathbb H^n$.

\paragraph{The case $e<0$ (that is, $\mathsf{K}=-2e/m>0$).}
The metric $\bar g_e$ is defined on all of $\mathbb R^n$. Stereographic
projection relates it to the metric $g_{_{\mathbb{S}^n}}$ on the sphere
$\mathbb{S}^n\subset\mathbb R^{n+1}$ centered at the origin and of curvature $\mathsf{K}$,
that is, of radius $R=1/\sqrt \mathsf{K}$. Let
$\phi:\mathbb{S}^n\setminus\{p\}\longrightarrow\mathbb R^n$ be stereographic
projection from the north pole $p=(0,\ldots,0,R)$ onto the equatorial
hyperplane $x_{n+1}=0$, and let $H_\mathsf{K}:\mathbb R^n\longrightarrow\mathbb R^n$
be the homothety of ratio $\mathsf{K}=1/R^2$. Then the composition\footnote{The
composition $\Phi$ agrees, in horizontal coordinates, with stereographic
projection onto the hyperplane $x_{n+1}=R-1/R$.}
$$
 \Phi:\mathbb{S}^n\setminus\{p\}\xrightarrow{\ \phi\ }\mathbb R^n
 \xrightarrow{\ H_\mathsf{K}\ }\mathbb R^n
$$
satisfies $\Phi_*g_{_{\mathbb{S}^n}}=\bar g_e$, where the pushforward of a metric
means $(\Phi^{-1})^*g_{_{\mathbb{S}^n}}$. Indeed, using the familiar identity
$$
 \phi_*g_{_{\mathbb{S}^n}}\,=\,\frac{4\sum_i \text{d}x_i^2}{(1+\mathsf{K}|x|^2)^2},
$$
we obtain
\begin{align*}
 \Phi_*g_{_{\mathbb{S}^n}}
 &\,=\,(H_\mathsf{K})_*(\phi_*g_{_{\mathbb{S}^n}})
 \,=\,(H_\mathsf{K})_*\left(\frac{4\sum_i \text{d}x_i^2}{(1+\mathsf{K}|x|^2)^2}\right)\\
 &\,=\,\frac{4\sum_i \text{d}x_i^2/\mathsf{K}^2}{(1+\mathsf{K}|x|^2/\mathsf{K}^2)^2}
 \,=\,\frac{4\sum_i \text{d}x_i^2}{(\mathsf{K}+|x|^2)^2}\,=\,\bar g_e\ .
\end{align*}
Applying $\Phi$ to the geodesics of $\mathbb{S}^n$, which are its great circles,
shows that \emph{the geodesics of $\bar g_e$ for $e<0$ are the circles and
straight lines in $\mathbb R^n$ that meet the sphere of radius
$\sqrt \mathsf{K}=1/R$ at antipodal points}.

\begin{center}
\begin{tikzpicture}[scale=1.25]
 \draw[step=.2,gray!15,very thin] (-2.25,-1.55) grid (1.65,2.2);
 \draw[densely dotted,thick] (0,0) circle (1.2);
 \draw[thick] (-.5,.5) circle ({sqrt(1.94)});
 \draw[thick] (-1.6,-1.6)--(1.65,1.65);
 \draw (0,0)--(1.2,0) node[midway,above] {$\sqrt \mathsf{K}$};
 \fill (0,0) circle (1.2pt) node[below left] {$0$};
 \fill ({1.2/sqrt(2)},{1.2/sqrt(2)}) circle (1.2pt);
 \fill ({-1.2/sqrt(2)},{-1.2/sqrt(2)}) circle (1.2pt);
 \node at (-1.0,1.6) {\small geodesic};
 \node at (1.08,1.92) {\small geodesic};
\end{tikzpicture}
\end{center}

\paragraph{The case $e=0$ (that is, $\mathsf{K}=0$).}
The metric $\bar g_e=4|x|^{-4}\sum_i \text{d}x_i^2$ is defined on
$\mathbb R^n\setminus\{0\}$. It is obtained from the usual Euclidean
metric $\sum_i \text{d}x_i^2$ by inversion with center the origin and radius
$\sqrt2$. \emph{The geodesics of $\bar g_e$ for $e=0$ are the circles and
straight lines through the origin}, with the origin itself omitted.

\paragraph{The case $e>0$ (that is, $\mathsf{K}=-2e/m<0$).}
The metric $\bar g_e$ is defined on the complement in $\mathbb R^n$ of
the closed ball $B$ given by $r\le\sqrt{|\mathsf{K}|}$. It is related to the
hyperbolic metric as follows. In Minkowski space $\mathbb R^{n+1}$,
with Lorentz metric $\text{d}x_0^2-\sum_{i=1}^n \text{d}x_i^2$, consider the
hyperboloid
$$
 \mathbb H^n\,=\,\left\{x_0^2-\sum_{i=1}^n x_i^2=R^2,\quad x_0>0\right\}\ ,
 \qquad R=\frac1{\sqrt{|\mathsf{K}|}}\ .
$$
Equip it with the Riemannian metric
$g_{_{\mathbb H^n}}=(\sum_{i=1}^n \text{d}x_i^2-\text{d}x_0^2)|_{\mathbb H^n}$, of constant
curvature $\mathsf{K}=-1/R^2$. Let
$\phi:\mathbb H^n\setminus\{p\}\longrightarrow\mathbb R^n$ be projection
from $p=(R,0,\ldots,0)$ onto $x_0=0$. Then the composition
$$
 \Phi:\mathbb H^n\setminus\{p\}\xrightarrow{\ \phi\ }\mathbb R^n
 \xrightarrow{\ H_\mathsf{K}\ }\mathbb R^n
$$
satisfies $\Phi_*g_{_{\mathbb H^n}}=\bar g_e$ on its image.\goodbreak

\emph{The geodesics of $\bar g_e$ for $e>0$ are the portions in
$\mathbb R^n\setminus B$ of circles and straight lines that meet the
sphere of radius $\sqrt{|\mathsf{K}|}=1/R$ orthogonally}. 
\begin{center}
\begin{tikzpicture}[scale=1.25]
 \draw[step=.2,gray!15,very thin] (-2.15,-1.85) grid (1.65,2.15);
 \draw[thick] (-.9,1.25) circle ({sqrt(.9*.9+1.25*1.25-1.2*1.2)});
 \draw[thick] (-1.85,-1.85)--(1.65,1.65);
 \fill[white] (0,0) circle (1.2);
 \draw[densely dotted,thick] (0,0) circle (1.2);
 \draw (0,0)--(1.2,0) node[midway,below] {$\sqrt{|\mathsf{K}|}$};
 \fill (0,0) circle (1.2pt) node[above left] {$0$};
 \node at (-.95,1.96) {\small geodesic};
 \node at (1.08,1.92) {\small geodesic};
 \node at (-1.34,-1.62) {\small geodesic};
\end{tikzpicture}
\end{center}

Theorem~5.1 and the preceding description of its geodesics yield
Hamilton's theorem: \emph{for the Kepler potential
$\mathsf{U}(r)=-mk/r$, hodographs are circular arcs, or intervals of
straight lines through the origin}.

\subsection*{Symmetries of the Kepler problem}

\refstepcounter{thm}\label{obs:spaceform-killing}
\noindent\textbf{\thethm \ \   Killing vectors on space-forms.}
Let $X=\mathbb{S}^n,\mathbb R^n$, or $\mathbb H^n$ be one of the three complete,
simply connected space-form models. The following facts are standard.
\begin{itemize}
 \item Once a point $p\in X$ is fixed, every Killing vector field $K$ on
 $X$ is determined by its first jet $j_p^1K$.
 \item The Lie algebra of Killing vector fields on $X$ has a basis of
 $\binom n2+n=\binom{n+1}2$ elements,
 $\{K_{ij},K_a\}_{1\le i<j\le n,\,1\le a\le n}$, whose first jets at $p$
 are
$$
  j_p^1K_{ij}\,=\,y_i\partial_{y_j}-y_j\partial_{y_i}\ ,
  \qquad j_p^1K_a\,=\,\partial_{y_a}\ ,
$$
 where $(y_1,\ldots,y_n)$ are normal coordinates centered at $p$.
 \item The first $\binom n2$ elements span the Lie subalgebra
 $\langle K_{ij}\rangle_{i<j}\simeq\mathfrak{so}(n)$ of Killing vector
 fields vanishing at $p$. This is the Lie algebra of the group $SO(n)$
 of orientation-preserving isometries of $X$ that fix $p$.

The last $n$ elements form a basis of a vector subspace
$\langle K_a\rangle_{a=1}^n$ of the Killing algebra. This subspace consists of the Killing vector fields $K$ satisfying $(\nabla K)_p=0$. It is not a Lie subalgebra,
except when $X=\mathbb R^n$.
\end{itemize}

\refstepcounter{thm}\label{explicit-killing}
\noindent\textbf{\thethm \ \  Killing vectors of the metrics $\bar g_e$.}
Consider the metrics from Theorem~5.1,
\[
 \bar g_e\,=\,\frac4{(|x|^2-2e/m)^2}\sum_i \text{d}x_i^2
 \,=\,\frac4{(|x|^2+\mathsf{K})^2}\sum_i \text{d}x_i^2\ ,
 \qquad |x|^2-\frac{2e}{m}>0\ ,
\]
where $\mathsf{K}=-2e/m$ is the curvature. As explained above, these metrics are
not complete, but become complete after adding a point $p$ at infinity.
Apply the observations of 5.3 with this point as the distinguished point.

The formula for $\bar g_e$ shows immediately that it is invariant under
$SO(n)$, acting by rotations about the origin of $\mathbb R^n$. Thus we
have the Lie subalgebra of Killing vector fields
$$
 \mathfrak{so}(n)\,=\,\langle K_{ij}^e\rangle_{1\le i<j\le n}\ ,
 \qquad K_{ij}^e\,=\,x_i\partial_{x_j}-x_j\partial_{x_i}\ .
$$
This $SO(n)$ is precisely the group of orientation-preserving isometries
of the completed metric fixing the point at infinity $p$, which explains
the notation $K_{ij}$ used in Observation~5.3.

Each coordinate axis $Ox_a$ is a geodesic of the completed metric passing
through $p$. Up to the choice of orientation, its unit tangent vector is
$$
 \frac{x_a^2+\mathsf{K}}{2}\partial_{x_a}.
$$
A routine calculation yields that this tangent field along $Ox_a$ extends to a
unique Killing vector field $K_a$ satisfying $(\nabla K_a)_p=0$, namely
$$
 K_a^e\,=\,\frac{\mathsf{K}-|x|^2}{2}\partial_{x_a}
       +x_a\sum_i x_i\partial_{x_i}.
$$

We therefore obtain the following result.
\begin{proposition}\label{prop:kepler-killing}
The Lie algebra of Killing vector fields of the metric
$$
 \bar g_e\,=\,\frac4{(|x|^2+\mathsf{K})^2}\sum_i \mathrm{d}x_i^2\ ,
 \qquad \mathsf{K}\,=\,-\frac{2e}{m}\ ,
$$
on $|x|^2+K>0$ has the basis
$$
 K_{ij}^e\,=\,x_i\partial_{x_j}-x_j\partial_{x_i}\ ,
 \qquad
 K_a^e\,=\,\frac{\mathsf{K}-|x|^2}{2}\partial_{x_a}
       +x_a\sum_i x_i\partial_{x_i}\ ,
$$
where $1\le i<j\le n$ and $1\le a\le n$.
\end{proposition}

\begin{remark} Let us give the Lie brackets of the preceding vectors (with $K_{ji}=-K_{ij}$, $K_{ii}=0$)
$$\begin{aligned}
[K_{ij},K_{k\ell}]
&\,=\,\delta_{jk}K_{i\ell}
-\delta_{ik}K_{j\ell}
-\delta_{j\ell}K_{ik}
+\delta_{i\ell}K_{jk},\\
[K_{ij},K_a]&\,=\,\delta_{ja}K_i-\delta_{ia}K_j,\\
[K_a,K_b]&\,=\,-\mathsf{K}\, K_{ab}\,=\,\frac{2e}{m}K_{ab}.
\end{aligned}
$$
\end{remark}

\refstepcounter{thm}
\noindent\textbf{\thethm\ \  Restricted symmetries of the Kepler problem.} By Theorem \ref{4.3}, the Killing vector fields  $\{K_{ij}^e,K_a^e\}$ of the metric $\bar g_e$ induce
 $H_e$-restricted symmetries $\{\mathbb{K}_{ij}^e,\mathbb{K}_a^e\}$ of the Kepler problem. Let us state without proof the values of these vector fields:
$$
\begin{aligned}
\mathbb{K}_{ij}^e
&\,=\,x_i\partial_{x_j}-x_j\partial_{x_i}
+v_i\partial_{v_j}-v_j\partial_{v_i}\ ,\\
\mathbb{K}_a^e
&\,=\,-\sum_i x_iv_i\partial_{x_a}
+\sum_i (x_av_i
-v_ax_i)\partial_{x_i}
-\left(\frac{e}{m}+\frac{|v|^2}{2}\right)\partial_{v_a}
+v_a\sum_i v_i\partial_{v_i}\ .
\end{aligned}
$$

\noindent\textit{Comments.}\smallskip

 (a). On each energy hypersurface
 $H_e\subset T(\mathbb R^n\setminus\{0\})$, the vector space
 $\langle \mathbb{K}_{ij}^e,\mathbb{K}_a^e\rangle_{i<j,\,a}$ is a Lie
 algebra isomorphic  to  $\langle K_{ij}^e,K_a^e\rangle_{i<j,\,a}$, hence isomorphic to $\mathfrak{so}(n+1)$ for $e<0$, to
 $\mathfrak{so}(1,n)$ for $e>0$, and to the Lie algebra
 $\mathfrak{se}(n)$ of Euclidean motions for $e=0$.\smallskip
 
  (b). The vector space $\langle \mathbb{K}_{ij}^e\rangle_{i<j}$ is a Lie
algebra isomorphic to $\mathfrak{so}(n)$, associated with the natural
 action of $SO(n)$ on $\mathbb R^n$. \smallskip
 
 (c). Let $\mathbb{K}_a$ be the vector field on $T(\mathbb{R}^n\setminus \{0\})$ defined by $(\mathbb{K}_a)_{|H_e}=\mathbb{K}_a^e$ for all $e\in\mathbb{R}$.  Although $\mathbb{K}_a^e$ is a $H_e$-restricted symmetry for each $e\in\mathbb{R}$, the vector field $\mathbb{K}_a$ is not a global symmetry on $T(\mathbb{R}^n\setminus \{0\})$, since $
L_{\mathbb K_a}\omega=\text{d}x_a\wedge \text{d}\textsf{E}$.
 On the contrary, the vector fields $\mathbb{K}_{ij}=x_i\partial_{x_j}-x_j\partial_{x_i}
+v_i\partial_{v_j}-v_j\partial_{v_i}$ are global symmetries.\smallskip

(d). The vector space of vector fields
 $\langle \mathbb{K}_{ij},\mathbb{K}_a\rangle_{i<j,\,a}$ on $T(\mathbb{R}^n\setminus \{0\})$  is not a Lie algebra.\bigskip

\refstepcounter{thm}
\noindent\textbf{\thethm\ \  Noether invariants of the Kepler problem.} Let us compute the Noether invariants on $H_e^*$ for the dual Kepler problem and then the Noether invariants on $H_e$ for the Kepler problem. By Theorem \ref{4.3}, the Killing vector fields  $\{K_{ij}^e,K_a^e\}$ of the metric $\bar g_e$ induce
 $H_e^*$-restricted symmetries $\{\mathbb{K}_{ij}^{e\,*},\mathbb{K}_a^{e\, *}\}$ of the  dual Kepler problem.
 Writing
$v_x=(x,v)\in T\mathbb R^n$, the corresponding Noether invariants on
$H_e^*$  are
$$
 h_{ij}^*(x,v)\,=\,h_{ij}^*(v_x)\,\overset{\ref{cor:jacobi-killing}}{=}\,g(v_x,K_{ij})
              \,=\,m(x_iv_j-x_jv_i),
$$
and
\begin{align*}
 h_a^*(x,v)&\,=\,h_a^*(v_x)\,\overset{\ref{cor:jacobi-killing}}{=}\,g(v_x,K_a)
       \,=\,m\frac{\mathsf{K}-|x|^2}{2}v_a+mx_a\sum_i x_iv_i\\
  (\text{substituting }\mathsf{K}=-2e/m)\\
 &\,=\,m\frac{-2e/m-|x|^2}{2}v_a+mx_a\sum_i x_iv_i\\
 &\,=\,-\left(e+\frac m2|x|^2\right)v_a+mx_a\sum_i x_iv_i\ .
\end{align*}
Applying the diffeomorphism
$J:T\mathbb R^n\longrightarrow T\mathbb R^n$, $(x,v)\mapsto(v,x)$,
to these functions gives Noether invariants, associated to the symmetries $\{\mathbb{K}_{ij}^e,\mathbb{K}_a^e\}$, on the energy hypersurface
$H_e$ of the original Kepler problem:
$$\boxed{\quad\begin{aligned}
 h_{ij}(x,v)&\,\overset{\ref{remark:dual-noether}}{=}\,-h_{ij}^*(v,x)\,=\,m(v_jx_i-v_ix_j),\\
 h_a(x,v)&\,\overset{\ref{remark:dual-noether}}{=}\,-h_a^*(v,x)
       \,=\,\left(e+\frac m2|v|^2\right)x_a-mv_a\sum_i x_iv_i\ .
\end{aligned}
}$$

Substituting $e=\mathsf{E}=\tfrac12m|v|^2-km/|x|$ into these $H_e$-restricted invariants gives
the following global invariants.

\begin{theorem}\label{thm:runge-lenz}
The following functions on $T(\mathbb R^n\setminus\{0\})$ are
invariants of the Kepler problem:
\begin{align*}
 h_{ij}(x,v)&\,=\,m(x_iv_j-x_jv_i)\ ,\\
 h_a(x,v)&\,=\,\left(m|v|^2-\frac{km}{|x|}\right)x_a
                  -mv_a\sum_i x_iv_i\ ,
\end{align*}
where $1\le i<j\le n$ and $1\le a\le n$.
\end{theorem}

\noindent\textbf{Note.} The invariants $h_{ij}$ are the components of the angular momentum tensor. The invariants $h_a$ are the components of the normalized Runge--Lenz vector 
$$
R(x,v):=\, m\left(|v|^2x-\langle x,v\rangle v-k\frac{x}{|x|}\right).
$$

The global invariants $h_{ij}$ are the Noether invariants associated to the global symmetries $\mathbb{K}_{ij}$ introduced in Subsection 5.6(c). As previously mentioned, the vector fields $\mathbb{K}_a$ are not global symmetries.

\end{document}